\documentstyle[12pt]{article}

\catcode`\@=11
\begin{document}
\centerline{\large\bf Drinfel$^{\prime}$d Realization of Quantum }
\centerline{\large\bf  Affine Superalgebra $U_q\widehat{(gl(1|1))}$}
\vspace{0.8cm}
\centerline{\sf $^a$Jin-fang Cai,  $^{bc}$S.K. Wang, $^a$Ke Wu 
  and $^a$Wei-zhong Zhao }
\baselineskip=13pt
\vspace{0.5cm}
\centerline{$^a$ Institute of Theoretical Physics, Academia Sinica, }
\baselineskip=12pt
\centerline{ Beijing, 100080, P. R. China }
\vspace{0.3cm}
%\centerline{ and}
%\vspace{0.3cm}
%\centerline{\sf Shi-kun Wang}
\baselineskip=13pt
\centerline{ $^b$ CCAST (World Laboratory), P.O. Box 3730, Beijing, 100080, 
  P. R. China  }
\baselineskip=12pt
\vspace{0.3cm}
\centerline{ $^c$ Institute of Applied Mathematics, Academia Sinica,}
\baselineskip=12pt
\centerline{Beijing, 100080, P. R. China }
\vspace{0.9cm}
\begin{abstract} {We obtain Drinfel$^{\prime}$d's realization of 
quantum affine superalgebra  $U_q\widehat{(gl(1|1))}$
 based on the super  version 
of RS construction method  and Gauss decomposition.}
\end{abstract}
\baselineskip=14pt
\vspace{1cm}

How to construct quantum algebras (including quantum affine algebras 
and Yangians) is a important problem in mathematics physics. After 
Drinfel$^{\prime}$d \cite{DRI1}
and Jimbo \cite{J1,J2}
 independently discovered that the universal enveloping 
algebra $U(g)$ of any simple Lie algebras or Kac-Moody
algebras admits a Hopf algebraic structure and a certain $q$-deformation, 
Drinfel$^{\prime}$d \cite{DRI2} gave his second definition or realization of 
quantum affine algebras $U_q(\widehat{g})$ and Yangians.
From views of quantum inverse scattering method, Faddeev, Reshetikhin 
and Takhtajan (FRT) \cite{FRT} gave another realization of $U_q(g)$ by 
means of a solution of Yang-Baxter equation (YBE). The FRT method has a much
more direct physics meaning and can be extended to quantum loop algebra
$U_q(g\otimes [t,t^{-1}])$ when the solution of YBE be dependent on  spectrum parameter. Laterly, Reshetikhin and 
Semenov-Tian-Shansky (RS)
 \cite{RS} used the exact affine anologue of FRT method to obtain
 a realization of quantum affine algebra $U_q(\widehat{g})$ . 
The explicit isomorphism between two realizations of quantum affine 
algebras $U_q(\widehat{g})$ given by Drinfel$^{\prime}$d and RS 
was established by Ding and Frenkel \cite{DF} using Gauss decomposition.
%-----------------------------------------------------------------
However, for a non-standard solution (without spectrum parameter) 
of YBE, if we employ the usual FRT constructive method, we would obtain 
a peculiar quantum algebra whose classical limit is not a Lie superalgebra
despite some of its relations (such as $ X^2=Y^2=0$) look like fermionic 
relations \cite{GE}. It was pointed out by Liao and Song \cite{LS} that,
to deal with a would-be quantum Lie superalgebra from  this nonstandard 
solution of YBE, one must start from formulas (YBE and RLL relations etc.) 
appropriate for the super case  at the very begining---super version of 
FRT method. We must use super version of FRT (RS) method to to construct 
quantum (affine) superalgebras \cite{LS, MR}.
Recently, some attentions have been paid to the construction of 
the  quantum affine superalgebras \cite{YA, FH}. 
In this paper, we will use the super affine version FRT 
 method (or super RS method) to construct quantum affine superalgebra
and use Gauss decomposition to get its Drinfel$^{\prime}$d realization. 
 We will focus on the simplest one: $U_q\widehat{(gl(1|1))}$, and this method 
can be easily extended to general case.

By $V$ we denote a ${\bf Z}_2$-graded two-dimension vector space 
(graded auxiliary space), and set the first and second basis of $V$ is 
even and odd respectively. In graded case, we must use this 
tensor product form: $(A \otimes B)(C \otimes D) =(-1)^{P(B)P(C)}
AC \otimes BD $, here $P(A)=0,1$ when $A$ is bosonic and Fermionic 
respectively. Then the graded(super) Yang-Baxter equation(YBE) takes
 this form\cite{LS,MR}: 

\begin{equation} 
\eta _{12}R_{12}(z/w) \eta _{13}R_{13}(z) \eta _{23}R_{23}(w)
= \eta _{23}R_{23}(w) \eta _{13}R_{13}(z) \eta _{12}R_{12}(z/w),
\end{equation}
where $R$-matrix acts in tensor product of two 2-D graded linear space $V$:
$R(z) \in End(V \otimes V) $, 
$\eta _{ik,jl}=(-1)^{P(i)P(k)} \delta _{ij} \delta _{lk}$ . The $R$-matrix 
must obey weight conservation condition: $R_{ij,kl}\neq 0 $
 only when $P(i)+P(j) = P(k)+P(l) $, here $P(1)=0$ and $P(2)=1$.
 It's very obvious that $\eta R(z)$ 
satisfies the ordinary YBE when $R(z)$ is a solution of super YBE. 
The super YBE can also be writed in components as:
\begin{eqnarray}
&&R^{ij}_{ab}(z/w)R^{ak}_{pc}(z)R^{bc}_{qr}(w)
(-1)^{(P(a)-P(p))P(b)}   \nonumber    \\
&&{\hspace{2cm}}=(-1)^{P(e)(P(f)-P(r))}
R^{jk}_{ef}(w)R^{if}_{dr}(z)R^{de}_{pq}(z/w)
\end{eqnarray} 

It's can be easily verified that following $R$-matrix  \cite{PS}
is a solution of graded YBE (1) :

\begin{equation} 
R_{12}(z)=\left(\begin{array}{lccr}
1&0 & 0& 0\\ 
0& \frac{z-1}{zq-q^{-1}} & \frac{z(q-q^{-1})}{zq-q^{-1}} & 0\\
0& \frac{(q-q^{-1})}{zq-q^{-1}} & \frac{z-1}{zq-q^{-1}} & 0 \\
0&0&0&-\frac{q-zq^{-1}}{zq-q^{-1}} 
\end{array} \right)
\end{equation}
This solution satisfy the unitary condition: $R_{12}(z)R_{21}(z^{-1})=
{\bf 1}$, and when $z=0$ and $q$ replaced by $q^{-1}$, the $\eta R(z)$
 degenerate to  the non-standard solution of YBE used in studying 
quantum superalgebra $U_q(gl(1|1))$\cite{LS, MR}. This solution $\eta R(z)$ 
can also be obtained  from the non-standard solution through Baxterization 
procedure \cite{YH}.

From the above solution  of graded YBE, we can define the quantum 
affine superalgebra  $U_q\widehat{(gl(1|1))}$ with
a central extension  employing super version of RS method or the 
affine version of what in\cite{LS}.  
 $U_q\widehat{(gl(1|1))}$ is an associative algebra with generators
\{$l_{ij}^k \vert 1\leq i,j \leq 2, k\in {\bf Z }$\} which subject to the 
following multiplication relations:
\begin{eqnarray}
&R_{12}(\frac{z}{w})L_1^{\pm}(z)  \eta L_2^{\pm}(w) \eta =
 \eta L_2^{\pm}(w) \eta L_1^{\pm}(z) R_{12}(\frac{z}{w})  
\\
&R_{12}(\frac{z_-}{w_+})L_1^{+}(z)  \eta L_2^{-}(w) \eta =
 \eta L_2^{-}(w) \eta L_1^{+}(z) R_{12}(\frac{z_+}{w_-})
\end{eqnarray}
here $z_{\pm}=zq^{\pm \frac{c}{2}}$. We have used standard notation:
$L_1^{\pm}(z)=L^{\pm}(z)\otimes {\bf 1}, L_2^{\pm}(z)={\bf 1} \otimes
L^{\pm}(z) $ and $L^{\pm}(z)=\left( l_{ij}^{\pm}(z) \right)^2_{i,j=1}$,
 $l_{ij}^{\pm}(z)$ are generating functions (or currents) of $l_{ij}^k$:
$l_{ij}^{\pm}(z)=\sum_{k=0}^{\infty}l_{ij}^{\pm k} z^{\pm k}$.

This algebra admits coalgebra and antipole structure compatible with the 
associative \\
 multiplication defined by eqs.(4) and (5): 

\begin{eqnarray}
&\triangle \left(l_{ij}^{\pm}(z)\right) =\sum_{k=1}^2 l_{kj}^{\pm}
(zq^{\pm \frac{c_2}{2}})\otimes 
l_{ik}^{\pm}(zq^{\mp \frac{c_1}{2}})(-1)^{(k+i)(k+j)}, \nonumber \\
& \epsilon\left(l_{ij}^{\pm}(z)\right)=\delta _{ij} ,{\hspace{1.cm}}
S\left( ^{st} L^{\pm}(z) \right) =\left[ ^{st} L^{\pm}(z) \right]^{-1}.
\end{eqnarray}
where  $\left[^{st} L^{\pm}(z) \right]_{ij} =(-1)^{i+j}l_{ji}^{\pm}(z) $.

 As what in Ding and Frenkel's paper \cite{DF},
 $L^{\pm}(z)$ have the following unique decompositions : 
\begin{eqnarray}
L^{\pm}(z)& = & \left( \begin{array}{lr} 1& 0\\f^{\pm}(z) & 1 \end{array}
\right) \left( \begin{array}{lr}k_1^{\pm}(z) & 0\\0&k_2^{\pm}(z) \end{array}
\right) \left( \begin{array}{lr} 1 & e^{\pm}(z) \\ 0&1 \end{array} 
\right)  \nonumber \\
&=&\left( \begin{array}{cc} k_1^{\pm}(z) & k_1^{\pm}(z)e^{\pm}(z) \\
f^{\pm}(z)k_1^{\pm}(z) & k_2^{\pm}(z) +f^{\pm}(z)k_1^{\pm}(z)e^{\pm}(z) 
\end{array} \right) 
\end{eqnarray}
where $ e^{\pm}(z), f^{\pm}(z) $ and $k_i^{\pm}(z)$ ($i$=1,2) is generating 
functions of  $U_q\widehat{(gl(1|1))}$ and $k_i^{\pm}(z)$ ($i$=1,2) are 
invertible. The inversions of $L^{\pm}(z)$  can be writed as:
\begin{equation}
L^{\pm}(z)^{-1}=\left( \begin{array}{cc}k_1^{\pm}(z)^{-1}+e^{\pm}(z)
k_2^{\pm}(z)^{-1}f^{\pm}(z) & -e^{\pm}(z)k_2^{\pm}(z)^{-1} \\
-k_2^{\pm}(z)^{-1}f^{\pm}(z) & k_2^{\pm}(z)^{-1} \end{array} \right)
\end{equation}
Let
\begin{equation}
X^+(z)=e^+(z_-)-e^-(z_+), {\hspace{1cm}} X^-(z)=f^+(z_+)-f^-(z_-).
\end{equation}

To calculate the (anti-)commutation relations of 
 $ X^{\pm}(z)$ and $k_i^{\pm}(z)$ ($i$=1,2) , we must 
make use of the the inversions of $L^{\pm}(z)$, unitarity of $R$-matirx
 and the following equivalent form of (4) and (5):
\begin{eqnarray}
&L_1^{\pm}(z)^{-1} \eta L_2^{\pm}(w)^{-1} \eta R_{12}(\frac{z}{w}) =
R_{12}(\frac{z}{w})  \eta L_2^{\pm}(w)^{-1} \eta L_1^{\pm}(z)^{-1}  \\
&L_1^{+}(z)^{-1} \eta L_2^{-}(w)^{-1} \eta R_{12}(\frac{z_-}{w_+}) =
R_{12}(\frac{z_+}{w_-})  \eta L_2^{-}(w)^{-1} \eta L_1^{+}(z)^{-1}  \\
&L_1^{\pm}(z) R_{12}(\frac{z}{w}) \eta L_2^{\pm}(w)^{-1} \eta = 
\eta L_2^{\pm}(w)^{-1} \eta  R_{12}(\frac{z}{w}) L_1^{\pm}(z)  \\
&L_1^{+}(z) R_{12}(\frac{z_+}{w_-}) \eta L_2^{-}(w)^{-1} \eta = 
 \eta L_2^{-}(w)^{-1} \eta  R_{12}(\frac{z_-}{w_+}) L_1^{+}(z)  \\
&L_1^{\pm}(z)^{-1} R_{21}(\frac{w}{z})  \eta L_2^{\pm}(w) \eta =
 \eta L_2^{\pm}(w) \eta R_{21}(\frac{w}{z}) L_1^{\pm}(z)^{-1}  \\
&L_1^{+}(z)^{-1} R_{21}(\frac{w_+}{z_-})  \eta L_2^{-}(w) \eta =
 \eta L_2^{-}(w) \eta R_{21}(\frac{w_-}{z_+}) L_1^{+}(z)^{-1}  
\end{eqnarray}

By the similar calculation process  made in \cite{DF} for quantum 
affine algebras , we can obtain the (anti-)commutation relations among 
$ X^{\pm}(z)$ and $k_i^{\pm}(z)$ ($i$=1,2).

From (4) (5) and (10)---(15), we can get all relations between 
$k_1^{\pm}(z)$ and $  k_2^{\pm}(z) $ :

\begin{eqnarray}
&& k_1^{\pm}(z) k_1^{\pm}(w)= k_1^{\pm}(w) k_1^{\pm}(z), {\hspace{1.0cm}}
  k_1^+(z) k_1^-(w)=k_1^-(w) k_1^+(z) \\
&& k_2^{\pm}(z)^{-1} k_2^{\pm}(w)^{-1}= k_2^{\pm}(w)^{-1} k_2^{\pm}(z)^{-1}  \\
&&\frac{w_+q-z_-q^{-1}}{z_-q-w_+q^{-1}} k_2^+(z)^{-1} k_2^-(w)^{-1} =
\frac{w_-q-z_+q^{-1}}{z_+q-w_-q^{-1}} k_2^-(w)^{-1} k_2^+(z)^{-1}   \\
&&k_1^{\pm}(z)k_2^{\pm}(w)=k_2^{\pm}(w)k_1^{\pm}(z)  \\
&&\frac{z_+-w_-}{z_+q-w_-q^{-1}}k_1^+(z)k_2^-(w)^{-1}=\frac{z_--w_+}
{z_-q-w_+q^{-1}}k_2^-(w)^{-1}k_1^+(z)   \\
&&\frac{w_+-z_-}{w_+q-z_-q^{-1}}k_2^+(z)^{-1}k_1^-(w)=\frac{w_--z_+}{w_-q-
z_+q^{-1}}k_1^-(w)k_2^+(z)^{-1}  
\end{eqnarray}

Then, we derive the relations between $k_1^{\pm}(z)$ and 
$X^{\pm}(w)$. From (4)(5) and unitarity of $R$-matrix, 
we have the following relations between $k_1^{\pm}(z)$ 
and $e^{\pm}(w), f^{\pm}(w)$:
\begin{eqnarray}
&&k_1^{\pm}(z)k_1^{\pm}(w)e^{\pm}(w)-\frac{z-w}{zq-wq^{-1}}k_1^{\pm}(w)
e^{\pm}(w)k_1^{\pm}(z)  \nonumber \\
&&{\hspace{5cm}}
-\frac{w(q-q^{-1})}{zq-wq^{-1}}k_1^{\pm}(w)k_1^{\pm}(z)e^{\pm}(z) =0 \\
&&k_1^{\pm}(z)k_1^{\mp}(w)e^{\mp}(w)-\frac{z_{\pm}-w_{\mp}}{z_{\pm}q-
w_{\mp}q^{-1}}k_1^{\mp}(w)e^{\mp}(w)k_1^{\pm}(z)  \nonumber \\
&&{\hspace{5cm}}-\frac{w_{\mp}(q-q^{-1})}{z_{\pm}q-w_{\mp}q^{-1}}k_1^{\mp}(w)
k_1^{\pm}(z)e^{\pm}(z) =0  \\
&&f^{\pm}(w)k_1^{\pm}(w)k_1^{\pm}(z)-\frac{z-w}{zq-wq^{-1}}k_1^{\pm}(z)
f^{\pm}(w)k_1^{\pm}(w) \nonumber \\
&&\hspace{5cm}-\frac{z(q-q^{-1})}{zq-wq^{-1}}f^{\pm}(z)k_1^{\pm}(z)k_1^{\pm}(w)=0 \\
&&f^{\mp}(w)k_1^{\mp}(w)k_1^{\pm}(z)-\frac{z_{\mp}-w_{\pm}}{z_{\mp}q-w_{\pm}
q^{-1}}k_1^{\pm}(z)f^{\mp}(w)k_1^{\mp}(w) \nonumber \\
&&\hspace{5cm}-\frac{z_{\mp}(q-q^{-1})}{z_{\mp}q-w_{\pm}q^{-1}}f^{\pm}(z)
k_1^{\pm}(z)k_1^{\mp}(w)=0 
\end{eqnarray}
thus
\begin{eqnarray}
&&(z_{\pm}-w)k_1^{\pm}(z)^{-1}e^{\pm}(w_{\mp})k_1^{\pm}(z)-
(z_{\pm}q-wq^{-1})e^{\pm}(w_{\mp})+w(q-q^{-1})e^{\pm}(z)=0 \\
&&(z_{\pm}-w)k_1^{\pm}(z)^{-1}e^{\mp}(w_{\pm})k_1^{\pm}(z)-(z_{\pm}q-wq^{-1})
e^{\mp}(w_{\pm})
+w(q-q^{-1})e^{\pm}(z)=0  \\
&&(z-w_{\pm})k_1^{\pm}(z)f^{\pm}(w_{\pm})k_1^{\pm}(z)^{-1}+z(q-q^{-1})
f^{\pm}(z)-(zq-w_{\pm}q^{-1})f^{\pm}(w_{\pm})=0   \\
&&(z_{\mp}-w)k_1^{\pm}(z)f^{\mp}(w_{\mp})k_1^{\pm}(z)^{-1}+z_{\mp}(q-q^{-1})
f^{\pm}(z)-(z_{\mp}q-wq^{-1})f^{\mp}(w_{\mp})=0  
\end{eqnarray}
So
\begin{eqnarray}
&&k_1^{\pm}(z)^{-1}X^+(w)k_1^{\pm}(z)=\frac{z_{\pm}q-wq^{-1}}
{z_{\pm}-w}X^+(w) \\
&&k_1^{\pm}(z)X^-(w)k_1^{\pm}(z)^{-1}=\frac{zq-w_{\pm}q^{-1}}
{z-w_{\pm}}X^-(w)
\end{eqnarray}

Now, we derive relations between  $k_2^{\pm}(z)$ and $X^{\pm}(w)$. 
From(10)(11) and unitarity of $R$-matirx, we have
\begin{eqnarray}
&&\frac{wq-zq^{-1}}{zq-wq^{-1}}e^{\pm}(z)k_2^{\pm}(z)^{-1}k_2^{\pm}(w)^{-1}
+\frac{z-w}{zq-wq^{-1}}k_2^{\pm}(w)^{-1}e^{\pm}(z)k_2^{\pm}(z)^{-1} \nonumber \\
&&{\hspace{6.5cm}}-\frac{z(q-q^{-1})}{zq-wq^{-1}}e^{\pm}(w)k_2^{\pm}(w)^{-1}
k_2^{\pm}(z)^{-1} =0\\
&&\frac{w_{\pm}q-z_{\mp}q^{-1}}{z_{\mp}q-w_{\pm}q^{-1}}e^{\pm}(z)k_2^{\pm}
(z)^{-1}k_2^{\mp}(w)^{-1}
+\frac{z_{\pm}-w_{\mp}}{z_{\pm}q-w_{\mp}q^{-1}}k_2^{\mp}(w)^{-1}e^{\pm}(z)
k_2^{\pm}(z)^{-1}\nonumber \\
&&{\hspace{6.5cm}}-\frac{z_{\pm}(q-q^{-1})}{z_{\pm}q-w_{\mp}q^{-1}}e^{\mp}(w)
k_2^{\mp}(w)^{-1}k_2^{\pm}(z)^{-1}=0  \\
&&\frac{z-w}{zq-wq^{-1}}k_2^{\pm}(z)^{-1}f^{\pm}(z)k_2^{\pm}(w)^{-1}-
\frac{w(q-q^{-1})}{zq-wq^{-1}}k_2^{\pm}(z)^{-1}k_2^{\pm}(w)^{-1}f^{\pm}(w)
\nonumber \\
&&\hspace{6.5cm}+\frac{wq-zq^{-1}}{zq-wq^{-1}}k_2^{\pm}(w)^{-1}
k_2^{\pm}(z)^{-1}f^{\pm}(z)=0 \\
&&\frac{z_{\mp}-w_{\pm}}{z_{\mp}q-w_{\pm}q^{-1}}k_2^{\pm}(z)^{-1}
f^{\pm}(z)k_2^{\mp}(w)^{-1}-
\frac{w_{\pm}(q-q^{-1})}{z_{\mp}q-w_{\pm}q^{-1}}k_2^{\pm}(z)^{-1}k_2^{\mp}
(w)^{-1}f^{\mp}(w)
\nonumber \\
&&\hspace{6.5cm}+\frac{w_{\mp}q-z_{\pm}q^{-1}}{z_{\pm}q-w_{\mp}q^{-1}}
k_2^{\mp}(w)^{-1}k_2^{\pm}(z)^{-1}f^{\pm}(z)=0 
\end{eqnarray}
Using the relations between $k_2^+(z)$ and $k_2^-(w)$ (17) and (18) ,
 then we get
\begin{eqnarray}
&&(w_{\pm}q-zq^{-1})e^{\pm}(z_{\mp})+(z-w_{\pm})k_2^{\pm}(w)^{-1}e^{\pm}
(z_{\mp})k_2^{\pm}(w)-z(q-q^{-1})e^{\pm}(w)=0  \\
&&(w_{\mp}q-zq^{-1})e^{\pm}(z_{\mp})+(z-w_{\mp})k_2^{\mp}(w)^{-1}e^{\pm}
(z_{\mp})k_2^{\mp}(w)-z(q-q^{-1})e^{\mp}(w)=0   \\ 
&&(z_{\pm}-w)k_2^{\pm}(w)f^{\pm}(z_{\pm})k_2^{\pm}(w)^{-1}-w(q-q^{-1})
f^{\pm}(w)+(wq-z_{\pm}q^{-1})f^{\pm}(z_{\pm})=0 \\
&&(z_{\mp}-w)k_2^{\mp}(w)f^{\pm}(z_{\pm})k_2^{\mp}(w)^{-1}-w(q-q^{-1})
f^{\mp}(w)+(wq-z_{\mp}q^{-1})f^{\pm}(z_{\pm})=0 
\end{eqnarray}
So
\begin{eqnarray}
&&k_2^{\pm}(z)^{-1}X^+(w)k_2^{\pm}(z)=\frac{z_{\pm}q-wq^{-1}}
{z_{\pm}-w}X^+(z)    \\
&&k_2^{\pm}(z)X^-(w)k_2^{\pm}(z)^{-1}=\frac{zq-w_{\pm}q^{-1}}
{z-w_{\pm}}X^-(w)
\end{eqnarray}

Now, we calculate the relations between  $X^{\pm}(z)$ 
and   $X^{\pm}(w)$ . From (4)(5) and unitarity of $R$-matirx , 
we have
\begin{eqnarray}
&&k_1^{\pm}(z)e^{\pm}(z)k_1^{\pm}(w)e^{\pm}(w)-\frac{wq-zq^{-1}}{zq-wq^{-1}}
k_1^{\pm}(w)e^{\pm}(w)k_1^{\pm}(z)e^{\pm}(z)=0 \\
&&k_1^+(z)e^+(z)k_1^-(w)e^-(w)-\frac{w_-q-z_+q^{-1}}{z_+q-w_-q^{-1}}
k_1^-(w)e^-(w)k_1^+(z)e^+(z)=0 \\
&&\frac{wq-zq^{-1}}{zq-wq^{-1}}f^{\pm}(z)k_1^{\pm}(z)f^{\pm}(w)k_1^{\pm}(w)-
f^{\pm}(w)k_1^{\pm}(w)f^{\pm}(z)k_1^{\pm}(z)=0  \\
&&f^+(z)k_1^+(z)f^-(w)k_1^-(w)+\frac{w_+q^{-1}-z_-q}{w_+q-z_-q^{-1}}
f^-(w)k_1^-(w)f^+(z)k_1^+(z)=0  
\end{eqnarray}
then, using the above relations between $k_1^{\pm}(z)$ and $X^{\pm}(w)$
(or $e^{\pm}(w)$ and $f^{\pm}(w)$), 
it's easily to get the following relations :
\begin{eqnarray}
&&X^+(z)X^+(w)+X^+(w)X^+(z)=0  \\
&&X^-(z)X^-(w)+X^-(w)X^-(z)=0
\end{eqnarray}

From (4)(5) and unitarity of $R$-matirx, we have following relations  :

\begin{eqnarray}
&&-(z-w)k_1^{\pm}(z)e^{\pm}(z)f^{\pm}(w)k_1^{\pm}(w)
-(z-w)f^{\pm}(w)k_1^{\pm}(w)k_1^{\pm}(z)e^{\pm}(z)\nonumber\\
&&{\hspace{4cm}}+z(q-q^{-1})\left(k_2^{\pm}(z)+f^{\pm}(z)k_1^{\pm}(z)e^{\pm}(z)
\right)k_1^{\pm}(w)  \nonumber  \\
&&{\hspace{4cm}}-z(q-q^{-1})\left(k_2^{\pm}(w)+f^{\pm}(w)k_1^{\pm}(w)e^{\pm}(w)
\right)k_1^{\pm}(z)=0  \\
&&\frac{z_{\mp}-w_{\pm}}{z_{\mp}q-w_{\pm}q^{-1}}k_1^{\pm}(z)e^{\pm}(z)
f^{\mp}(w)k_1^{\mp}(w)
+\frac{z_{\pm}-w_{\mp}}{z_{\pm}q-w_{\mp}q^{-1}}f^{\mp}(w)k_1^{\mp}(w)
k_1^{\pm}(z)e^{\pm}(z) \nonumber \\
&&{\hspace{4cm}}-\frac{z_{\mp}(q-q^{-1})}{z_{\mp}q-w_{\pm}q^{-1}}
\left(f^{\pm}(z)k_1^{\pm}(z)e^{\pm}(z)+k_2^{\pm}(z)
\right)k_1^{\mp}(w)  \nonumber \\
&&{\hspace{4cm}}+\frac{z_{\pm}(q-q^{-1})}{z_{\pm}q-w_{\mp}q^{-1}}
\left(f^{\mp}(w)k_1^{\mp}(w)e^{\mp}(w)+k_2^{\mp}(w)
\right)k_1^{\pm}(z) =0 
\end{eqnarray}
then, using the following relations also from(4)(5) and 
unitarity of $R$-matirx:
\begin{eqnarray}
&&k_1^{\pm}(z)e^{\pm}(z)k_1^{\pm}(w)-\frac{z(q-q^{-1})}{zq-wq^{-1}}
k_1^{\pm}(w)e^{\pm}(w)k_1^{\pm}(z) \nonumber \\
&&{\hspace{6cm}}-\frac{z-w}{zq-wq^{-1}}k_1^{\pm}(w)k_1^{\pm}
(z)e^{\pm}(z)=0  \\
&&\frac{z-w}{zq-wq^{-1}}k_1^{\pm}(z)f^{\pm}(w)k_1^{\pm}(w)
-f^{\pm}(w)k_1^{\pm}(w)k_1^{\pm}(z)  \nonumber  \\
&&{\hspace{6cm}}+\frac{z(q-q^{-1})}{zq-wq^{-1}}f^{\pm}(z)
k_1^{\pm}(z)k_1^{\pm}(w)=0  
\end{eqnarray}
and
\begin{eqnarray}
&&k_1^{\pm}(z)e^{\pm}(z)k_1^{\mp}(w)-\frac{z_{\pm}(q-q^{-1})}
{z_{\pm}q-w_{\mp}q^{-1}}k_1^{\mp}(w)e^{\mp}(w)k_1^{\pm}(z) \nonumber \\
&&{\hspace{5cm}}-\frac{z_{\pm}-w_{\mp}}{z_{\pm}q-w_{\mp}q^{-1}}
k_1^{\mp}(w)k_1^{\pm}(z)
e^{\pm}(z)=0  \\
&&\frac{z_{\mp}-w_{\pm}}{z_{\mp}q-w_{\pm}q^{-1}}k_1^{\pm}(z)f^{\mp}(w)
k_1^{\mp}(w)-f^{\mp}(w)k_1^{\mp}(w)k_1^{\pm}(z) \nonumber \\
&&{\hspace{5cm}}+\frac{z_{\mp}(q-q^{-1})}{z_{\mp}q-w_{\pm}q^{-1}}
f^{\pm}(z)k_1^{\pm}(z)k_1^{\mp}(w)
=0
\end{eqnarray}
we obtain:
\begin{eqnarray}
&&(z-w)\left(e^{\pm}(z)f^{\pm}(w)+f^{\pm}(w)e^{\pm}(z)\right) \nonumber \\
&&{\hspace{2.5cm}}=z(q-q^{-1})\left(k_1^{\pm}(z)^{-1}k_2^{\pm}(z)
-k_1^{\pm}(w)^{-1}k_2^{\pm}(w)\right) \\
&&(z_{\mp}-w_{\pm})\left(e^{\pm}(z)f^{\mp}(w)+f^{\mp}(w)e^{\pm}(z)\right)
\nonumber  \\
&&{\hspace{2.5cm}}=z_{\mp}(q-q^{-1})k_1^{\pm}(z)^{-1}k_2^{\pm}(z)-
\frac{z_{\mp}-w_{\pm}}{z_{\pm}-w_{\mp}}z_{\pm}(q-q^{-1})
k_1^{\mp}(w)^{-1}k_2^{\mp}(w)
\end{eqnarray}
we must mind that the expansion direction of (54) can be chosed in
 $\frac{w}{z}$ or $\frac{z}{w}$, but for the first and second of (55), 
the expansion direction is only in $\frac{w}{z}$  and $\frac{z}{w}$
respectively. Then we can get:

\begin{eqnarray}
&&X^+(z) X^-(w)+X^-(w)X^+(z)  \nonumber \\
&=&(q-q^{-1})\left
[\delta\left(\frac{w_-}{z_+}\right)k_1^-(z_+)^{-1}k_2^-(z_+)-
\delta\left(\frac{z_-}{w_+}\right)k_1^+(w_+)^{-1}k_2^+(w_+)\right]
\end{eqnarray}
here $\delta (z)=\sum_{k\in {\bf Z}}z^k$.

In a word, we get  all relations between $k_i^{\pm}(z)$
 (i=1,2) and $X^{\pm}(w)$  as follows: 

\begin{eqnarray}
&&[k_1^{\pm}(z) ~,~ k_1^{\pm}(w)]=[k_1^+(z) ~,~ k_1^-(w)]=0, \\
&&[k_1^{\pm}(z) ~,~ k_2^{\pm}(w)]=[k_2^{\pm}(z) ~,~ k_2^{\pm}(w)]=0 , \\
&&\frac{w_+ q-z_- q^{-1}}{z_- q - w_+q^{-1}}k_2^+(z)^{-1}k_2^-(w)^{-1}=
  \frac{w_- q-z_+ q^{-1}}{z_+ q - w_-q^{-1}}k_2^-(w)^{-1}k_2^-(z)^{-1} , \\
&&\frac{z_{\pm}-w_{\mp}}{z_{\pm}q-w_{\mp}q^{-1}}k_1 ^{\pm}(z) k_2^{\mp}(w)^{-1} =
\frac{z_{\mp}-w_{\pm}}{z_{\mp}q-w_{\pm}q^{-1}}k_2 ^{\mp}(w)^{-1} k_1^{\pm}(z) ,  \\
&&k_i^{\pm}(z)^{-1}X^+(w)k_i^{\pm}(z)=\frac{z_{\pm}q-wq^{-1}}{z_{\pm}-w}X^{+}(w) ,
{\hspace{0.6cm}} (i=1,2) \\
&&k_i^{\pm}(z)X^-(w)k_i^{\pm}(z)^{-1}=\frac{z_{\mp}q-wq^{-1}}{z_{\mp}-w}X^{-}(w) ,
{\hspace{0.6cm}} (i=1,2) \\
&&\{ X^+(z) ~,~ X^+(w) \}=\{ X^-(z) ~,~ X^-(w) \}=0  \\
&&\{ X^+(z) ~,~ X^-(w) \}=(q-q^{-1})\left[ \delta(\frac{w_-}{z_+})
k_1^-(z_+)^{-1}k_2^-(z_+)-\delta(\frac{z_-}{w_+})k_1^+(w_+)^{-1}k_2^+(w_+) \right] 
\end{eqnarray}

The above relations are  Drinfel$^{\prime}$d's construction of 
quantum affine superalgebra $U_q\widehat{(gl(1|1))}$. It's very clear that 
 $ X^{\pm}(z)$ (or $ e^{\pm}(z)$ and $f^{\pm}(z) $ ) are Fermionic type for
their anti-commutation relations and 
$k_i^{\pm}(z)$ ($i$=1,2) are Bosonic type elements in  $U_q\widehat{(gl(1|1))}$
as expected. This result differs from what in \cite{FH} for the relation 
between $ X^+(z)$ and $ X^-(w)$ (64) is also anti-commutator which is 
requirement of superalgebra.

Introducing a transformation for the generating functions
$ X^{\pm}(z)$ and $k_i^{\pm}(z)$ ($i$=1,2) :
\begin{eqnarray}
&&E(z)=X^+(zq) {\hspace{2.5cm}} F(z)=X^-(zq) \\
&&K^{\pm}(z)=k_1^{\pm}(zq)^{-1}k_2^{\pm}(zq) {\hspace{0.7cm}} 
H^{\pm}(z)=k_2^{\pm}(zq)k_1^{\pm}(zq^{-1})
\end{eqnarray}
then the (anti-)commutation relations for $E(z), F(z), K^{\pm}(z)$
and $H^{\pm}(z) $ are
\begin{eqnarray}
&&[K^{\pm}(z) ~,~K^{\pm}(w)]=[H^{\pm}(z) ~,~ H^{\pm}(w)] =0 \\
&&[K^+(z) ~,~ K^-(w)] =[K^{\pm}(z) ~,~H^{\pm}(z)]=0 \\
&&H^+(z)H^-(w)=\left(\frac{(z_-q-w_+q^{-1})(z_+q^{-1}-w_-q)}
{(z_+q-w_-q^{-1})(z_-q^{-1}-w_+q)}\right)^2 H^-(w)H^+(z)  \\
&&K^{\pm}(z)H^{\mp}(w)=\frac{(w_{\mp}q-z_{\pm}q^{-1})(z_{\mp}q-w_{\pm}q^{-1})}
{(w_{\pm}q-z_{\mp}q^{-1})(z_{\pm}q-w_{\mp}q^{-1})} H^{\mp}(w)K^{\pm}(z)  \\
&&[K^{\pm}(z) ~,~ E(w)]=[K^{\pm}(z) ~,~ F(w)]=0  \\
&&E(w)H^{\pm}(z)=\frac{z_{\pm}q-wq^{-1}}{z_{\pm}q^{-1}-wq} H^{\pm}(z)E(w) \\
&&F(w)H^{\pm}(z)=\frac{z_{\mp}q^{-1}-wq}{z_{\mp}q-wq^{-1}} H^{\pm}(z)F(w) \\
&&\{E(z) ~,~E(w) \}=\{F(z) ~,~F(w) \}=0  \\
&&\{ E(z) ~,~ F(w) \}=(q-q^{-1})\left[\delta\left(\frac{w_-}{z_+}\right)
K^-(z_+)-\delta\left(\frac{z_-}{w_+}\right)K^+(w_+) \right] 
\end{eqnarray}

Set $e^{\pm}(z_{\mp}q)=E^{\pm}(z)$ and $f^{\pm}(z_{\pm}q)=F^{\pm}(z)$, 
then $ E(z)=E^+(z)-E^-(z)$ and $F(z)=F^+(z)-F^-(z) $. 
The comultiplication of the generating functions $E^{\pm}(z),
F^{\pm}(z), K^{\pm}(z)$ and  $H^{\pm}(z)$ can be calculated as 
follows:
\begin{eqnarray}
\triangle \left(K^{\pm}(z)\right)&= &K^{\pm}(zq^{\pm \frac{c_2}{2}}) \otimes 
 K^{\pm}(zq^{\mp \frac{c_1}{2}})  \\
\triangle \left(H^{\pm}(z)\right) &= &H^{\pm}(zq^{\pm \frac{c_2}{2} }) \otimes 
 H^{\pm}(zq^{\mp \frac{c_1}{2} }) -q^{-1}F^{\pm}(zq^{\mp \frac{c_1}{2}\pm 
\frac{c_2}{2}}) H^{\pm}(zq^{\pm \frac{c_2}{2}})\otimes  \\
&& \left(E^{\pm}(zq^{\mp \frac{c_1}{2}\pm \frac{c_2}{2}})
H^{\pm}(zq^{\mp\frac{c_1}{2}})+H^{\pm}(zq^{\mp\frac{c_1}{2}})E^{\pm}
(zq^{\mp \frac{c_1}{2}\pm \frac{c_2}{2}}) \right)  \\
\triangle \left(E^{\pm}(z)\right)&=&E^{\pm}(z) \otimes
 {\bf 1}+K^{\pm}(zq^{\mp\frac{c_1}{2}})\otimes E^{\pm}(zq^{\mp c_1})  \\
\triangle \left(F^{\pm}(z)\right)&=&{\bf 1}\otimes F^{\pm}(z) 
+F^{\pm}(zq^{\pm c_2})\otimes K^{\pm}(zq^{\pm\frac{c_2}{2}}) 
\end{eqnarray}

 The general quantum affine superalgebras $U_q\widehat{(gl(m|n))}$ 
will be discussed in another paper and the analogous
procedure can also be applied in studying super Yangian doubles using a 
rational solution of super YBE \cite{C2}.

{\bf Acknowledge:} The authors would like to thank Dr. X.M. Ding 
for informing us the Ref.\cite{FH} before publication.

\end{document}